\documentclass[12pt,a4paper]{article}
\usepackage{latexsym,amsmath,amssymb,graphics,amscd}
\usepackage{graphics,color}

\usepackage{array}

\newtheorem{theorem}{Theorem}[section]

\usepackage{latexsym,amsmath,amssymb,graphics,amscd}
\usepackage{graphicx}

\newcommand{\EE}{\mathbb{E}}
\newcommand{\var}{\mathsf{var}}

\newcommand{\E}{\mathbb{E}}

\newcommand{\PP}{\mathbb{P}}

\newcommand{\RR}{\mathbb{R}}
\newcommand{\hhf}{{\scriptstyle{{\frac{1}{2}}}}}

\title{An Introduction to Multilevel Monte Carlo for Option Valuation\thanks{Submitted to 
International Journal of
Computer Mathematics, special issue on  Computational Methods in Finance} } 

\author{
Desmond J. Higham\thanks{Department of Mathematics and Statistics, University of Strathclyde, UK} 
}

\usepackage{array}

\begin{document}
\maketitle

\begin{abstract}
Monte Carlo is a simple and flexible tool that is widely used in computational 
finance.
In this context, it is common for the 
quantity of interest to be the expected value of a random variable defined via a stochastic differential equation. 
In 2008, Giles proposed a remarkable  
improvement to the 
approach of discretizing with a numerical method and applying 
standard Monte Carlo. His 
\emph{multilevel Monte Carlo} method 
offers a speed up of $O(\epsilon^{-1})$, where $\epsilon$ is the required accuracy.
So computations can run 100 times more quickly 
when two digits of accuracy are required.
The \lq\lq multilevel philosophy\rq\rq\ has since been adopted by a range of 
researchers and a wealth of practically significant
results has arisen, most of which have yet to make their way into the expository
literature.
  In this work, we give a brief, accessible, introduction to multilevel Monte Carlo and summarize
recent  results 
 applicable to the task of option evaluation.
\end{abstract}


\section{Aims}
\label{sec:aims}

Finding the appropriate market value of a financial option can usually be formulated as 
an expected value computation \cite{Glasserman,Hbook,Sey12}. 
In the case where the product underlying the option is modelled as a 
stochastic differential equation (SDE), we may  
\begin{itemize}
 \item simulate the SDE numerically to compute many independent sample paths, and then 
 \item combine the option payoff from each path in order to obtain a Monte Carlo estimate, and 
 an accompanying confidence interval.
\end{itemize}
Compared with other approaches, notably the direct discretization of a partial differential equation based formulation of the problem, a Monte Carlo computation has the advantages of
 (a) being simple to implement and (b) being flexible enough to cope with a wide range of   
underlying SDE models and option payoffs. On the downside, 
Monte Carlo is typically expensive in terms of computation time \cite{Glasserman,Hbook}.

In the seminal 2008 paper \cite{Giles08}, Giles pulled together ideas from numerical analysis,
stochastic analysis and 
applied statistics in order to deliver a dramatic improvement on the efficiency of the
\lq\lq
SDE simulation plus Monte Carlo\rq\rq\ approach.
If the required level of accuracy, in terms of confidence interval, is $\epsilon$, the multilevel
approach essentially improves the computational complexity by a factor of $\epsilon$. So
for a calculation requiring two digits of accuracy, we obtain a hundredfold improvement in  
computation time.   
Multilevel Monte Carlo has rapidly become an extremely hot topic in the field of 
 stochastic computation, impacting on a wide range of application areas. In particular, technical reviews of research progress in the field have 
begun to appear 
\cite{Giles14,GS12}
and 
a comprehensive survey is currently in progress by Giles for the journal Acta Numerica.
However, the area is still sufficiently new that most
textbooks in computational finance do not introduce the topic, and hence it has 
not been fully integrated into typical graduate-level classes and development courses for practitioners.
For this reason, we present here 
an accessible introduction to the multilevel Monte Carlo
 approach, and give a brief overview of the current state of 
the art with respect to financial option valuation.    

In section~\ref{sec:conv} we discuss the underlying SDE simulation.
 Section~\ref{sec:mc} then considers the complexity of standard Monte Carlo in this setting.
 In section~\ref{sec:mlmot}
 we give some motivation for the multilevel approach, which is introduced and analysed in 
section~\ref{sec:emmlmc}.
Section~\ref{sec:compex} illustrates the performance of the algorithm in practice, using 
code that has been made available by Giles.
In section~\ref{sec:follow} we give pointers to multilevel research in option valuation 
that has built on \cite{Giles08}.
Section~\ref{sec:disc} concludes with a brief discussion.

\section{Convergence in SDE Simulation}
\label{sec:conv}

We consider an Ito SDE
of the form
\begin{equation}
 dX(t) = f(X(t)) dt + g(X(t)) dW(t), \quad X(0) = X_0.
 \label{eq:sde}
 \end{equation}
Here, $f:\RR^m \to \RR^m$ and 
$g:\RR^m \to \RR^{m \times d}$
are given functions, known as the drift and diffusion coefficients,
respectively, 
and $W(t) \in \RR^{d}$ is standard Brownian motion.
The initial condition $X_0$ is supplied and we wish to simulate the SDE over the fixed time interval $[0,T]$.
The Euler--Maruyama method \cite{KP99,MT04} 
computes approximations 
$X_n \approx X(t_n)$, where $t_n = n \Delta t$, 
 according to  $X_0 = X(0)$ and, for $n  = 1,2,\ldots N-1$, 
\begin{equation}
 X_{n+1} = X_n +  f(X_n)  \Delta t + g(X_n) \Delta W_n,
\label{eq:em}
\end{equation}
where $\Delta t = T/N$ is the stepsize and $\Delta W_n = W(t_{n+1}) - W(t_n)$ is the relevant  
Brownian motion increment.

In the study of the accuracy of SDE simulation methods, the  
two most widely used convergence concepts are referred 
to as \emph{weak} and \emph{strong}
\cite{KP99,MT04}.
Roughly, 
\begin{itemize}
 \item weak convergence controls the error of the means, whereas,
 \item strong convergence controls the mean of the error.
\end{itemize}
To prove weak and strong convergence results, 
we must impose conditions on the SDE. For example it is standard 
to assume that $f$ and $g$ in (\ref{eq:sde}) satisfy global Lipschitz 
conditions; that is, there exists a constant $L$ such that 
\begin{equation}
| f(x) - f(y) | \le L | x - y| 
\quad 
\mathrm{and}
\quad 
| g(x) - g(y) | \le L | x - y|, 
\quad 
\mathrm{for~all~} x,y \in \RR^m.
\label{eq:gl}
\end{equation}
Here and throughout we take $\| \cdot \|$ to be the Euclidean norm.
Under such conditions, and for appropriate initial data, it follows that  the 
Euler--Maruyama method has weak order one, so that
\begin{equation}
 \sup_{0 \le t_n  \le T}
   \left(
   \EE [ X(t_n) ] -
   \EE [ X_n ]
     \right)
       = O( \Delta t).
   \label{eq:emweak}
   \end{equation}
   In the sense of strong error, which involves the mean of the absolute difference between
   the two random variables at each grid point, Euler--Maruyama achieves only an order of one half in general:
   \begin{equation}
   \EE \left[
 \sup_{0 \le t_n \le T}
    \left|
       X(t_n) -
    X_n \right| \right]
       = O( \Delta t^{\hhf} ).
   \label{eq:strong}
   \end{equation}
More generally, for any $m > 1$ and sufficiently small $\Delta t$ 
there is a constant $C = C(m)$ such that 
   \begin{equation}
   \EE \left[
 \sup_{0 \le t_n \le T}
    \left|
       X(t_n) -
    X_n \right|^m \right]
   \le C \Delta t^{m/2}.
   \label{eq:strongm}
   \end{equation}
Strong convergence is sometimes described as a pathwise property.
This can be understood via the Borel-Cantelli Lemma. 
 For example,
 in \cite{KN07} it is shown that
given any $\epsilon > 0$ there exists a path-dependent constant
$K = K(\epsilon)$ such that, for all sufficiently small $\Delta  t$,
\[
 \sup_{0 \le t_n \le T}
    \left|
       X(t_n) -
    X_n \right|
     \le K (\epsilon) h^{\hhf - \epsilon}.
   \]

In the setting of this work, it is tempting to argue that strong convergence is not relevant; 
if
we wish to compute an expected value
based on the SDE solution then following individual paths accurately is not important.
However, we will see in section~\ref{sec:emmlmc}
that the analysis in   
\cite{Giles08} justifying multilevel Monte Carlo 
makes use of both weak and strong convergence properties.

To conclude this section, we remark that the analysis of SDE simulation on problems that 
violate the global Lipschitz conditions (\ref{eq:gl}) is far from complete.
In the case of SDE models for financial assets and interest rates,
issues may arise through faster than linear growth
at infinity and also through unbounded derivatives at the origin.
For example, both complications occur in the class of scalar interest rate models from 
  \cite{A99},
\[
dX(t) = \left(\alpha_{-1} X(t)^{-1} - \alpha_0 + \alpha_1 X(t) - 
       \alpha_2 X(t)^r \right) dt + \alpha_3 X(t)^\rho dW(t),
\]
where the $\alpha_i $ are positive constants and $r, \rho > 1$.
        Although some positive results are available for specific nonlinear structures 
         \cite{HMS02,HJ11,HJ14,SMHP10},
 there has also been a sequence of negative results
showing how Euler--Maruyama
  can break down
on nonlinear SDEs 
\cite{HMS02,HJ09,MT05}.

\section{SDE Simulation and Standard Monte Carlo}
\label{sec:mc}

Given the SDE
(\ref{eq:sde}),
suppose we wish to approximate the final time
 expected value of the solution, $\EE[X(T)]$, using Monte Carlo
 with Euler--Maruyama. 
 We will let $\epsilon$ denote the accuracy requirement in terms of confidence interval width; fixing on a 95\%\  confidence level to be concrete, we therefore wish to be in a position where 
applying the algorithm independently a large number of times, the exact answer would be within
$\pm \epsilon$ of our computed answer with frequency at least $0.95$.

Let $X_N^{[s]}$ denote the
 Euler--Maruyama final time approximation along the $s$th path.
Using $M$ Monte Carlo samples 
we may form the sample average
 \[
  a_M = \frac{1}{M} \sum_{s=1}^{M} X_N^{[s]}.
  \]
 The overall error in our approximation has the form 
 \begin{eqnarray}
   a_M - \EE[X(T)] &=&
              a_M - \EE[ X(T) - X_N + X_N]      \nonumber \\
                &=& a_M - \EE[X_N] + \EE[X_N - X(T) ]. \label{eq:mcsplit}
              \end{eqnarray}
   Note that $X_N$ denotes a random variable describing the
   result of applying Euler--Maruyama
 (\ref{eq:em}), whereas each $X_N^{[s]}$ is an independent sample from
 the distribution given by $X_N$.  
    The expression  (\ref{eq:mcsplit}) breaks down the error into two terms.
   The statistical error, $a_M - \EE[X_n]$,  is concerned with how well
 we can approximate an expected value from a finite number of samples;  
           it does not
               depend on how accurately the numerical method approximates the SDE (in
                  particular it does not depend significantly on $\Delta t$) and it will
                   generally decrease if we take more sample paths.
 The discretization error, or bias, 
$\EE[X_N - X(T) ]$,
  arises because we have 
 approximated the SDE with a difference equation; this  
is the discrepancy that would remain if we had access to the exact
    expected value of the numerical solution and it will generally decrease
      if we reduce the stepsize.
    
Standard results  
\cite{Glasserman,RC04}
tell us that  
 the 
  statistical error $a_M - \EE[X_n]$
 can be described via
 a confidence interval of width $O(1/\sqrt{M})$.
 The weak convergence property (\ref{eq:emweak})
 shows that the bias 
$\EE[X_N - X(T) ]$ 
behaves like $O(\Delta t)$; so we must 
        add this amount to the confidence interval width.
        We arrive at an overall confidence interval of width
        $O(1/\sqrt{M}) + O(\Delta t)$. To achieve our required target accuracy
        of $\epsilon$, we see that 
   $ 1/\sqrt{M}$ and $\Delta t$
  should scale like $\epsilon$. In other words,
         $M$ should scale like $\epsilon^{-2}$ and $\Delta t$ should scale like
          $\epsilon$.

 It is reasonable to measure computational  cost by counting either the number of times that the drift and diffusion coefficients, $f$ and $g$,
        are evaluated, or
       the number of times that a random number generator is called. In either case,
 the cost per path is proportional to $1/\Delta t$, and hence the 
  total cost of the computation scales like 
  $M / \Delta t$.
  We argued above that  $M$ should scale like $\epsilon^{-2}$ and $\Delta t$ should scale like
          $\epsilon$.  Here is the conclusion:
 \begin{quote}
 we may achieve accuracy $\epsilon$ 
 by combining Euler--Maruyama and standard Monte Carlo 
 at an overall cost that 
 scales like $\epsilon^{-3}$.
  \end{quote}

One approach to improving the computational complexity is to 
replace Euler--Maruyama with a simulation method of higher weak order 
\cite{AM10,KP99,MT04}.
If we use a second order method, so that 
(\ref{eq:emweak}) is replaced by 
\[
 \sup_{0 \le t_n  \le T}
   \left(
   \EE [ X(t_n) ] -
   \EE [ X_n ]
     \right)
       = O( \Delta t^2),
\]
 then a straightforward adaption of the arguments above 
 lead to the following conclusion:
  \begin{quote}
 we may achieve accuracy $\epsilon$ 
 by combining a second order weak method and standard Monte Carlo 
 at an overall cost that 
 scales like $\epsilon^{-2.5}$.
  \end{quote}
We note, however, that establishing second order weak convergence requires extra 
smoothness assumptions to be placed on the SDE coefficients.

As we show in section~\ref{sec:emmlmc}, 
the method of Giles \cite{Giles08} has the following 
feature:
 \begin{quote}
 we may achieve accuracy $\epsilon$ 
 by using Euler--Maruyama 
   in a multilevel Monte Carlo setting 
 at an overall cost that 
 scales like $\epsilon^{-2} (\log \epsilon)^2 $.
  \end{quote}
Moreover, by using a higher strong order method, such as Milstein \cite{KP99,MT04}, it 
is possible to reduce the multilevel Monte Carlo cost to 
 the order of $\epsilon^{-2}$ \cite{Giles07}.

It is worth pausing to admire 
an $O(\epsilon^{-2})$ computational complexity 
count.
Suppose we are given an exact expression for the SDE solution, as a function of $W(t)$.
Hence, we are able to compute exact samples, without the need to apply a numerical method.
A standard Monte Carlo approach requires $1/\sqrt{M}$  to scale like $\epsilon$ in order
to achieve the required confidence interval width. 
If we regard the evaluation of each exact $X(T)$ sample as having $O(1)$ cost, 
then the cost overall will be proportional to $M$; that is, $\epsilon^{-2}$. 
In this sense, with a multilevel approach \emph{the numerical analysis comes for free}; 
we can solve the problem as quickly as one for which we 
have an exact pathwise expression for the SDE solution.

\section{Motivation for the Multilevel Approach}
\label{sec:mlmot}
\begin{figure}[htb]
\begin{center}
 \scalebox{0.7}{\includegraphics{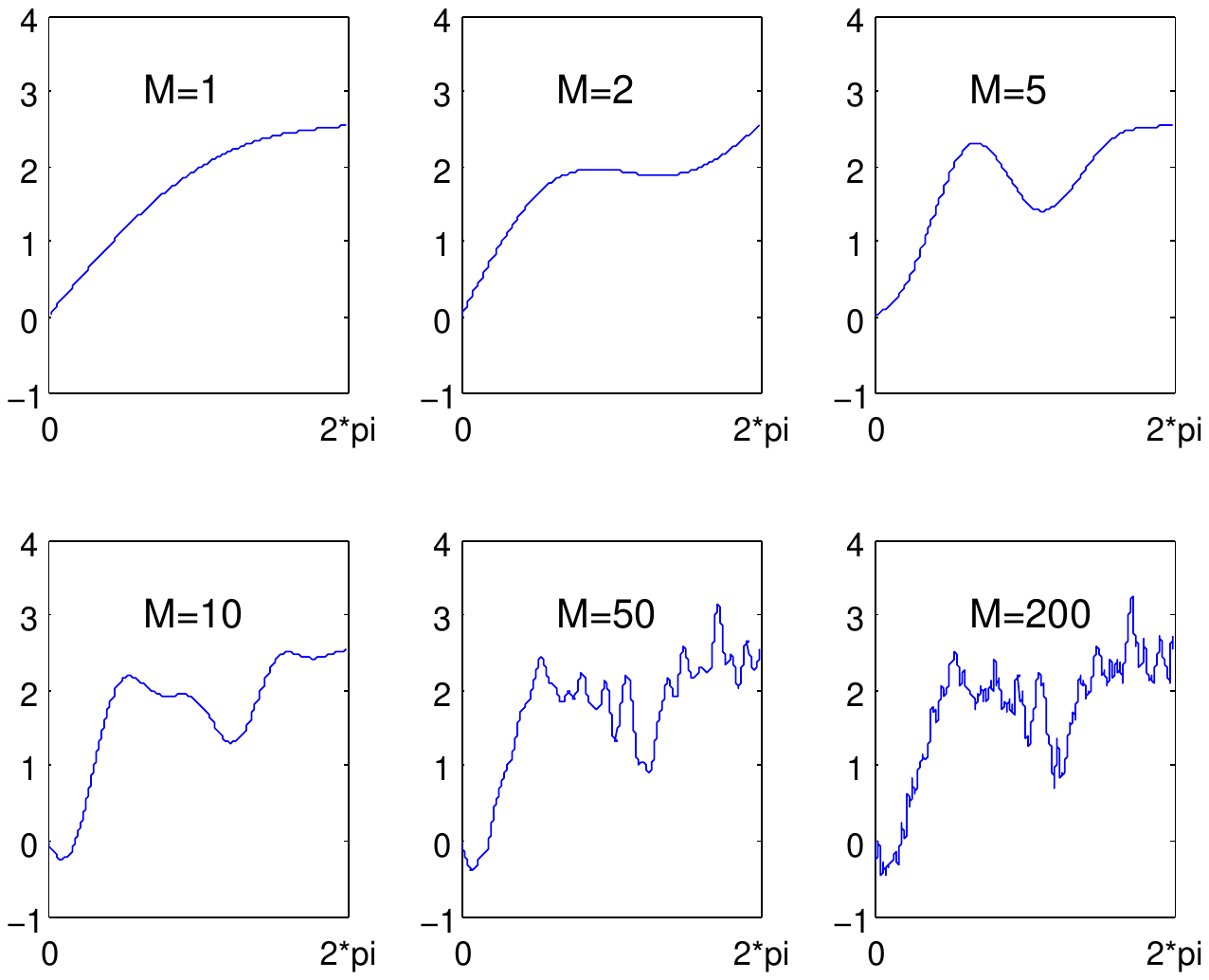}}
\end{center}
        \caption[bif]
                  {
    Paths based on
   the 
  Paley-Wiener representation
(\ref{eq:pw}).
  As indicated, the six plots show the
  sum truncated after
  $M = 1$, $2$,
  $5$, $10$,
  $50$ and  $200$ sine terms.
                \label{Fig.pw}
                }
\end{figure}

We can motivate the multilevel approach by considering 
a series expansion of Brownian motion, where the
coefficients are random variables.
The \emph{Paley-Wiener representation} 
over 
     $[0,2 \pi]$
has the form
\begin{equation}
   W(t) = Z_0  \frac{t}{\sqrt{2 \pi}}
            +
              \frac{2}{\sqrt{\pi}} \sum_{n=1}^{\infty}
               Z_n \frac{ \sin( \hhf n t) }{n},
 \label{eq:pw}
\end{equation}
where the $\{ Z_i \}_{i \ge 0}$ are i.i.d.\ and $N(0,1)$; see, for example, 
\cite{Mikosch}.
In Figure~\ref{Fig.pw}
we
draw samples for the $Z_i$ and plot the curves arising
when the infinite series in
(\ref{eq:pw}) is truncated to $\sum_{n=1}^{M}$, for
$M = 1$, $2$,
$5$, $10$,
$50$ and $200$.
It is clear that 
the early terms in the series
affect the overall shape, while the later terms add fine detail.
From this perspective, it is intuitively reasonable that 
we can build up information at different resolution scales, with the finer scales 
having less impact on the overall picture.

Now, we may view Monte Carlo as requiring a \lq\lq black box\rq\rq\  that
returns independent samples. In our numerical SDE context, 
the samples come from a distribution that is only approximately correct, and
the black box (the Euler--Maruyama method) comes with a dial. 
Turning the dial corresponds to changing $\Delta t$. 
Samples with a smaller $\Delta t$ are more expensive---we have to wait longer for them because
the paths contain more steps. The multilevel Monte Carlo algorithm 
cleverly exploits this dial. The black box is used to produce samples across a
range of stepsizes. Most of the samples that 
we ask for will be obtained quickly with relatively large $\Delta t$ values.
Relatively few 
samples will be generated at the expensive small $\Delta t$ levels.
In a sense, the large $\Delta t$ paths cover the low-frequency information so that
expensive, high-frequency paths are used sparingly. Figure~\ref{Fig.pw}
might convince you that this idea has some merit.
The next section works through the details.

\section{Multilevel Monte Carlo with Euler--Maruyama}
\label{sec:emmlmc}

We focus now on the more general case where we wish to approximate the expected value of
some function of the final time solution, $\EE[h(X(T))]$.
We have in mind the case where 
$X(t)$ represents an underlying asset, in risk-neutral form, and 
$h(\cdot)$ is the payoff of a corresponding 
European-style option
\cite{Glasserman,Hbook}.
For example, 
$h(x) = \max(x - E,0)$ for a European call option with exercise price $E$ and expiry time $T$.
For simplicity, we will consider the  scalar case, so that
$m = d = 1$ in 
 (\ref{eq:sde}), but we note that all arguments generalise to the case of systems, with the same 
conclusions.
We assume that the payoff function $h$ satisfies a global Lipschitz
condition; this covers the call and put option cases.

Multilevel Monte Carlo uses a range of different 
discretization levels.
At level $l$ we have a
stepsize of the form
\begin{equation}
 \Delta t_l = M^{-l} T, \quad \mathrm{where~} l = 0,1,2,\ldots,L.
 \label{eq:mlmc:Dtl}
 \end{equation}
Here $M>1$ is a fixed quantity whose precise value does not affect the overall
complexity of the method, in terms of the asymptotic rate as $\epsilon \to 0$.
For simplicity we may think of $M = 2$.
As the upper limit on the level index we choose 
 \begin{equation}
  L = \frac{\log \epsilon^{-1}}{\log M}.
   \label{eq:mlmc:Ldef}
  \end{equation}
  In this way, at the coarsest level, $l = 0$, we have the largest stepsize, 
    $\Delta t_0 = T$, covering the whole interval in one step.
 At the most refined level, $ l = L$, we have 
 $\Delta t_L = O(\epsilon)$---from (\ref{eq:emweak}), this the stepsize needed by 
 Euler--Maruyama to achieve weak error of $O(\epsilon)$.

\begin{figure}[htb]
\begin{center}
 \scalebox{0.7}{\includegraphics{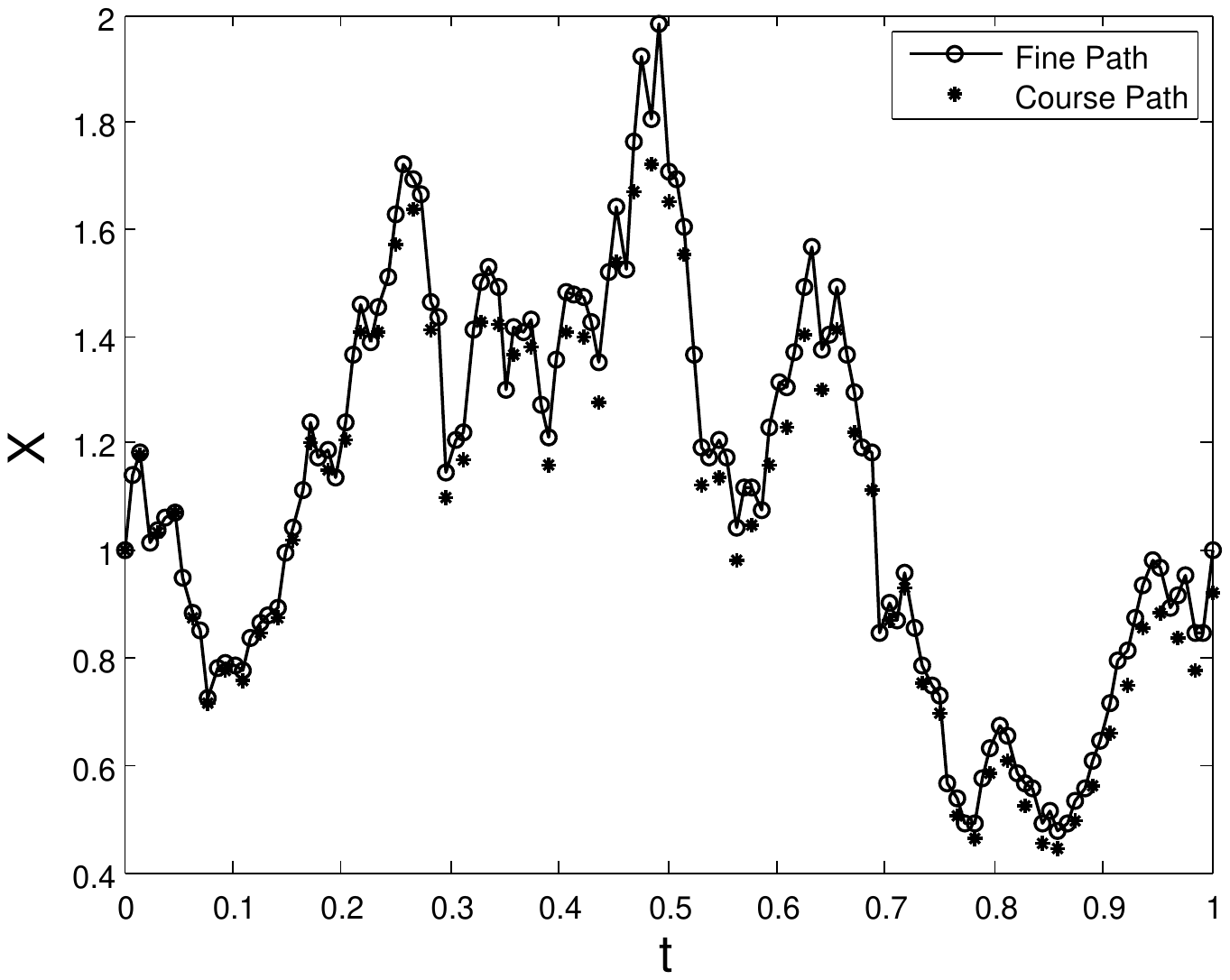}}
\end{center}
        \caption[bif]
                  {
 Illustration of how the estimator 
        $\widehat Y_l$
 in 
                        (\ref{eq:Yldef}) is constructed.
 Circles (joined by straight lines for clarity) 
 show the refined Euler--Maruyama path, with stepsize 
 $\Delta t = 2^{-l} T$.
 Asterisks 
 show the coarser Euler--Maruyama path, with stepsize 
 $\Delta t = 2^{-l+1} T$,  \emph{computed with the same Brownian increments}.
                \label{Fig.pair}
                }
\end{figure}

 With each choice of stepsize, $\Delta t_l$, we may apply Euler--Maruyama to the SDE
  (\ref{eq:sde}) and evaluate the payoff function $h$ at the final time.
  We will let the random variable $\widehat P_l$ denote this approximation to $h(X(T))$.
  Now, from the linearity of the expectation operator we have 
 the telescoping sum 
  \begin{equation}
      \EE[\widehat P_L] = \EE[\widehat P_0] + \sum_{l=1}^L \EE[\widehat P_l - \widehat P_{l-1}].
      \label{eq:mlmc:ident}
      \end{equation}
 In multilevel Monte Carlo, we 
       use the expansion on the right hand side as an indirect means to evaluate the 
   left hand side.
  This may be thought of as a recursive application of the \emph{control variate} technique, 
  which is widely used in applied statistics \cite{Glasserman,Hbook,RC04,Sey12}. 
  To estimate $\EE[\widehat P_0]$
  we form the usual sample mean, based on, say, $N_0$, paths. This gives 
   \begin{equation}
       \widehat Y_0 = \frac{1}{N_0} \sum_{s=1}^{N_0} \widehat P_0^{[s]}.
       \label{eq:mlmc:Yzero}
    \end{equation}
       Generally, for  $\EE[\widehat P_l - \widehat P_{l-1}]$ with $l > 0$ 
  we will use $N_l$ paths
  so that 
    \begin{equation}
        \widehat Y_l = \frac{1}{N_l} \sum_{s=1}^{N_l} \left(
                        \widehat P_l^{[s]} - \widehat P_{l-1}^{[s]} \right).
                        \label{eq:Yldef}
  \end{equation}
It is vital to point out that 
  $  \widehat P_l^{[s]} $ and $ \widehat P_{l-1}^{[s]}$ 
 in (\ref{eq:Yldef}) 
come from the \emph{same discretized Brownian path},
 with different stepsizes $\Delta t_l$ and $\Delta t_{l-1}$, respectively.
 Figure~\ref{Fig.pair} illustrates the idea for the case $M=2$.
 In words, at a general level $l$, we compute
$N_l$ Brownian paths and, for each path, apply Euler--Maruyama twice; once with stepsize
$\Delta t_l$ and 
once with stepsize $\Delta t_{l-1}$. (In practice, 
  we compute a path at resolution $\Delta t_l$ and then combine Brownian increments
  over pairs of steps in order to get a path at resolution $\Delta t_{l-1}$.)
  Having constructed our $N_l$ independent paths for level $l$, we start afresh at level $l+1$;
  none of the earlier information is re-used and new (independent) pseudo-random numbers
  are generated. 

  Because of the choice of $L$ in 
 (\ref{eq:mlmc:Ldef})
    we know from (\ref{eq:mlmc:ident}) that our estimator will have the required 
   $O(\epsilon)$ bias. Now we will see how to choose the values of $\{ N_l \}_{l=0}^{L}$
   to achieve the corresponding accuracy in the overall confidence interval.

 Considering a general level where $l > 0$, 
 appealing to the strong convergence behaviour
 (\ref{eq:strongm}) 
 of Euler--Maruyama 
 and our assumption that
 $h$ is globally Lipschitz, we have
\begin{eqnarray}
  \var[ \widehat P_l - h(X(T)) ] &=&
                           \EE [ ( \widehat P_l - h(X(T)) )^2 ]
                                  - (\EE [  \widehat P_l - h(X(T)) ] )^2 \label{a1} \\
    &\le& \EE [ ( \widehat P_l - h(X(T)) )^2 ] \label{a2} \\
    &\le& \mathrm{constant} \times  \EE [ (X_N - X(T) )^2 ]  \label{a3} \\
   &=& O(\Delta t_l). \label{eq:mlmc:varPlbd}
\end{eqnarray}
It then follows from Minkowski's Inequality \cite{CK}
 that
\begin{eqnarray}
  \var[ \widehat P_l - \widehat P_{l-1} ] &=& \var[ \widehat P_l - h(X(T)) + h(X(T))  - \widehat P_{l-1} ] 
                     \nonumber \\
              &\le& \left( \sqrt{ \var[ \widehat P_l - h(X(T))] } 
                + \sqrt{ \var[ h(X(T))  - \widehat P_{l-1} ] } \right)^2 \nonumber \\
      &=&  O(\Delta t_l).
\label{eq:mlmc:varPlm1}
\end{eqnarray}
Applying this result in 
(\ref{eq:Yldef}) 
we conclude that 
$\widehat Y_l$ 
has a variance of $O(\Delta t_l/N_l)$ for
$l > 1$. Because all levels are independent, we deduce 
that
\[
 \var\left[
        \widehat Y_0
      + 
    \sum_{l = 1}^{L}  
        \widehat Y_l
\right] = \var[ \widehat Y_0] + \sum_{l=1}^{L} O(\Delta t_l/N_l).
\]
To balance the variance evenly across levels $l = 1,2,\ldots,L$ and to 
control the variance at level $l = 0$, we 
choose $N_l = O(\epsilon^{-2} L \Delta t_l)$. It then follows that
our overall estimator has variance 
\[
 O(\epsilon^{2}) + \sum_{l=1}^{L} O(\epsilon^{2}/L)
                     = O(\epsilon^{2}).
\]
In this way, we have achieved the bias and variance required to give a confidence 
interval of the specified $\epsilon$ level of accuracy.

To quantify how the complexity of this algorithm scales with $\epsilon$,
we sum the cost of level $l$ from 
 $l = 0$ to $L$  
to give
\[
 \sum_{l=0}^{L} N_l \Delta t_l^{-1} =
 \sum_{l=0}^{L} \epsilon^{-2} L  \Delta t_l \Delta t_l^{-1} = L^2  \epsilon^{-2}.
\]
From (\ref{eq:mlmc:Ldef}) this  expression becomes 
$O(\epsilon^{-2} \left( \log \epsilon \right)^2)$, as we quoted in 
section~\ref{sec:mc}.

At this stage, a few remarks are in order:
\begin{description}
\item[Constructive Upper Bound:]
In the course of the analysis above,
we came up with a general-purpose choice for
 the
number of paths at each level, $\{ N_l \}_{l = 0}^{L}$.
The final complexity count is therefore an upper bound on the
best possible value. In practice, for a given problem and accuracy requirement,
we can perform
a cheap pre-processing step where appropriate variances are estimated and
an optimization problem is solved in order to give a
sequence  $\{ N_l \}_{l = 0}^{L}$; see, for example, 
\cite{GS12}. 

\item[Weak versus Strong:]
The key inequality 
(\ref{eq:mlmc:varPlm1}),
which guarantees tight coupling between coarse and fine paths,
made use of the strong convergence property. For small $\Delta t_l$, both paths are close to the
true path,
so the paths must be close to each other.
In this sense, both strong and weak error rates are key ingredients in the analysis.
We note, however, that 
Giles
\cite{Giles07} has also developed estimators that do not rely directly on strong convergence.

\item[Variance and Second Moment]
In deriving the  inequality
(\ref{eq:mlmc:varPlbd}), 
we discarded the square of the first moment and used the
second moment as an upper bound for the variance.
This may appear to be a very crude step, but
in our context it does not
degrade the final conclusion.
(In a different, Poisson-driven setting where a multilevel
 method was developed and analysed,
 the step  
(\ref{a1})--(\ref{a2}) 
is no longer optimal---it
is beneficial to analyse the variance directly 
\cite{AHS14}.)

\item[Exploiting Structure:]
As mentioned above,
multilevel Monte Carlo may be viewed as a recursive version of
the control variate approach. 
In the simplest version of control variates,  if we wish to compute
$\EE[ X ]$, we may instead
 compute $\EE[ X - Y]$ and add $\EE[Y]$, where $Y$ is a suitably constructed random variable
such that  $X - Y$ has small variance and $\EE [Y]$ is readily available
\cite{Glasserman,Hbook}. 
However, the success of this technique usually relies on
incorporating some extra knowledge of the problem: 
a structure such as symmetry or convexity, or the existence of
a \lq\lq nearby\rq\rq\  problem that is analytically tractable.
In this respect, the multilevel Monte Carlo method for SDEs is very different from
traditional control variates: the analysis is completely general and no special
insights are needed about the nature of the underlying SDE, other than
knowledge of the basic weak and strong
convergence properties.

\item[Multilevel versus Multigrid:]
In \cite{Giles14}, Giles explains that the \emph{multigrid} approach in numerical PDEs was
\lq\lq the inspiration for the author in developing the MLMC method for SDE path simulation.\rq\rq\ 
There are clear similarities between the two: the use of geometrically refined/coarsened grids 
and the idea that
relatively little work needs to  be expended on the fine grids in order 
to resolve high frequency components.
However, it is important to keep in mind that 
there are also conceptual differences between the two techniques: 
multilevel Monte Carlo is distinct, and novel. For example, multilevel Monte
Carlo does not involve the notion of passing information
up and down the refinement levels, as is done with multigrid V or W cycles.

\item[Related Earlier Methods:]
As discussed in \cite[section~1.3]{Giles14},
 related earlier work on improving Monte Carlo when samples are generated via discretization
was performed by Heinrich, see, for example
\cite{Hein1,Hein2}, 
and
Kebaier
\cite{Keb05}
devised a two-level approach to path simulation.

\end{description}

Based on the type of analysis that we summarized above, it is 
possible to state a general theorem about multilevel simulation:
\begin{theorem}[Giles; see for example, \cite{Giles14}]
Let $P$ denote a random variable, and let $P_l$ denote the corresponding 
level $l$ numerical approximation. If there exist  independent estimators $Y_l$ based
 on $N_l$ Monte Carlo samples, and positive constants 
 $\alpha, \beta, \gamma, c_1, c_2, c_3$ such that $\alpha \ge \hhf \min (\beta, \gamma)$  and
\begin{enumerate}
 \item  $| \EE[ P_l - P ] | \le c_1 2^{-\alpha l}$
 \item $  \EE[ Y_0] =  \E[ P_0]$ and \ $  \EE[ Y_l] =  \E[ P_l - P_{l-1}]$ for $l > 0$
 \item $\var[Y_l] \le c_2 N_l^{-1} 2^{-\beta l}$
 \item $  \EE[ C_l] \le c_3 N_l 2^{\gamma l} $, where $C_l$ is the computational complexity of $Y_l$
\end{enumerate}
then there exists 
 a positive constant $c_4$ such that for any $\epsilon < e^{-1}$ there are values $L$
and $N_l$ for which the  multilevel estimator 
\[
 Y = \sum_{l = 0}^{L} Y_l
\]
has a mean-square error with bound 
\[
 \EE\left[ (Y - \EE[P])^2\right] < \epsilon^2
\]
 with a computational complexity $C$ with bound 
\[
 \EE[ C ]  \le
 \left\{
 \begin{array}{ll}
  c_4 \epsilon^{-2},   & \beta > \gamma\\
  c_4 \epsilon^{-2} (\log(\epsilon))^2,   & \beta = \gamma\\
  c_4 \epsilon^{-2 - (\gamma - \beta)/\alpha},   & \beta < \gamma.
 \end{array}
 \right.
\]
\end{theorem}

Giles 
\cite{Giles07}
has also shown how to construct estimators for which $\beta > \gamma = 1$, by replacing
Euler--Maruyama with the more strongly accurate Milstein scheme. For European-style options with
Lipschitz payoff functions, this makes $O(\epsilon^{-2})$ complexity achievable.
From the arguments in section 
(\ref{sec:mc}), it is intuitively reasonable that this is the optimal rate.
The issue is formalized in 
\cite{mr09},
and optimality is confirmed.

In section~\ref{sec:conv}
 we mentioned that the basic Euler--Maruyama method 
(\ref{eq:sde})
 may fail to converge in a weak or strong sense on 
nonlinear SDEs in the asymptotic limit $\Delta t \to 0$.
A closely related question, of direct relevance to this review, is whether  the combination of
\lq\lq Euler--Maruyama plus Monte Carlo\rq\rq\ converges in the $\epsilon \to 0$ limit.
In 
\cite{HJ11,HJ14},
Hutzenthaler and Jentzen showed that Euler--Maruyama Monte Carlo
can achieve convergence in a $\PP$-almost sure sense in cases where the underlying
Euler--Maruyama scheme diverges.
  This can happen when the events causing Euler--Maruyama to diverge are so rare
 that they are extremely unlikely to impact on any of the Monte Carlo samples.
 However, 
in 
\cite{HJK13}
 Hutzenthaler, Jentzen
and Platen
showed that the multilevel Monte Carlo method does not inherit this property. 
They established this result using a
counterexample of the form
\begin{equation}
  dX(t) = -X(t)^5 dt,
\label{eq:count}
\end{equation}
with $X(0)$ having a standard Gaussian distribution, where $\EE[ X(t)^2 ]$ is the required moment.
	Note that the SDE
	(\ref{eq:count})
	has a zero drift term, so it may also be regarded as a random ODE.
	A modified version of Euler--Maruyama, known as a \emph{tamed} method,
	was shown in
	\cite{HJK13}
	to recover convergence in the multilevel setting.

	\section{Computational Experiments}
	\label{sec:compex}

	Asymptotic, $\epsilon \to 0$, analysis 
	indicates that multilevel Monte Carlo offers a dramatic improvement 
	in computational complexity. 
	Numerous computational studies have confirmed that this potential can be realised in 
	practice.

\begin{figure}[htb]
\begin{center}
 \scalebox{0.3}{\includegraphics{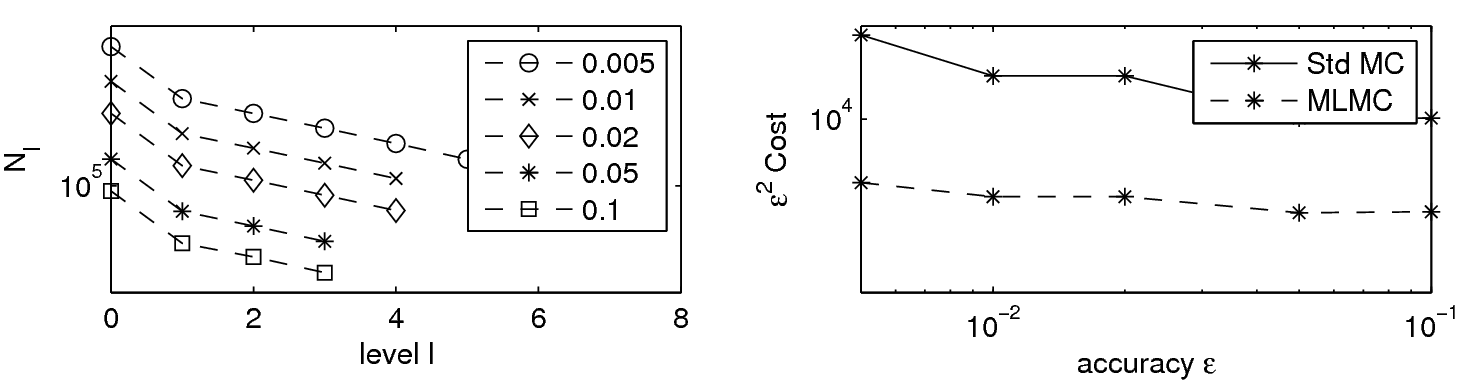}}
 \scalebox{0.3}{\includegraphics{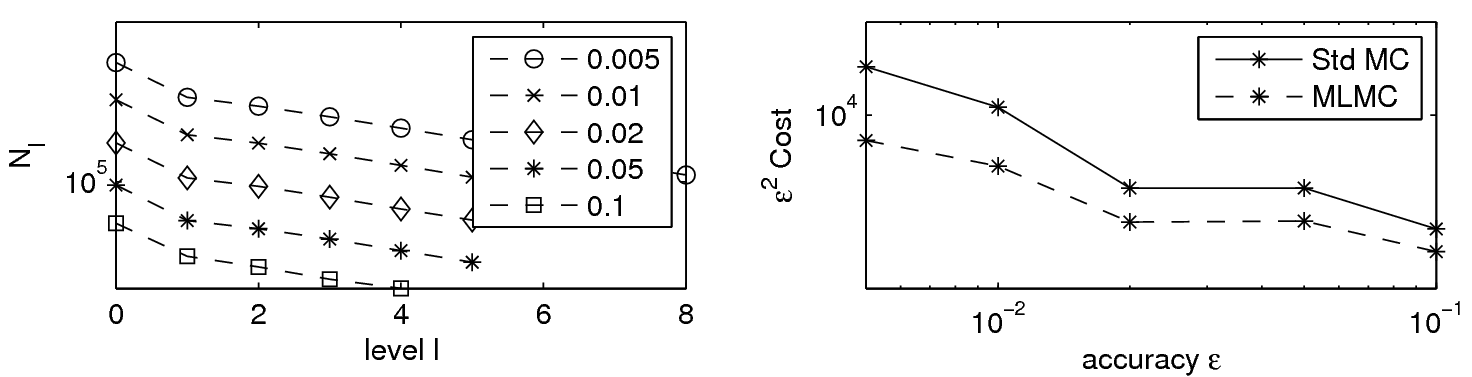}}
\end{center}
        \caption[bif]
                  {
 Output from the multilevel Monte Carlo code made available by Giles (see text for 
 web site address).
 Left hand pictures show the number of paths per level at each target accuracy.
 Right hand pictures show the computation time, scaled by $\epsilon^2$.
 Upper pictures are for a European call option.  
 Upper pictures are for a digital option.  
                \label{Fig.gilescode}
                }
\end{figure}

	Giles has made MATLAB code available at 
	\begin{verbatim}
	http://people.maths.ox.ac.uk/gilesm/acta/
	\end{verbatim}
	that can be used as the basis for computational experimentation.
	In Figure~\ref{Fig.gilescode} we show results 
	based on this code.
	Here, we have an asset model given by geometric Brownian motion
\[
 dX(t) =  0.05 X(t) dt + 0.25 X(t)  dW(t), \quad X(0) = 100.
\]
We consider (a) a European call and (b) a digital  
call option over $[0,T]$ with $T = 1$ and  exercise price 
$100$. So the payoff functions, after discounting for interest, are
\[
h(x) = e^{-0.05T} \max(x - 100,0)
\]
 for the call option and 
\[
h(x) = \left\{ \begin{array}{ll}
                       e^{-0.05T} 100 & \mathrm{~when~} x > 100 \\
                       0 & \mathrm{~when~} x < 100 
                 \end{array}
        \right.
\]
for the digital option.
(For those who worry about probability zero events, 
the code defines $h(100)  =  e^{-0.05T} (100 + 0)/2$.)
The code repeats the Monte Carlo simulation for accuracy requests of 
$\epsilon = 0.1, 0.05, 0.02, 0.01, 0.005$.
The upper left picture in 
	Figure~\ref{Fig.gilescode} shows, for the call option, 
the number of paths 
 $N_l$ used at each level $l$  in the multilevel method.
We see that for a given $\epsilon$  more paths are used at the cheaper 
(small $l$) levels, and as $\epsilon$ is decreased, so that more accuracy 
is required, extra levels are added.
The upper right picture indicates the corresponding computational cost in terms of
run time. 
More precisely,  the asterisks (joined by dashed lines) show the 
cost weighted by $\epsilon^2$ as a function of $\epsilon$.
We see that this quantity remains approximately constant, as predicted by the analysis.
The picture also shows the scaled cost for an equivalent 
standard Monte Carlo computation, using a solid linetype.
We see a much larger cost that appears to grow faster than $\epsilon^{-2}$.
The lower pictures in 
	Figure~\ref{Fig.gilescode}  give the same results for the case of the digital 
option, and again the multilevel version is seen to be more efficient than 
standard Monte Carlo.

\section{Follow-on Research}
\label{sec:follow}
In this section we summarize some of the key advances that have been
made since the original multilevel breakthrough \cite{Giles08}.
We focus on work that is directly relevant to financial option valuation.
The comprehensive overviews 
\cite{Giles14,GS12} 
can be consulted for
further details on these, and other, areas. 
The webpage maintained by Giles at
\begin{verbatim}
http://people.maths.ox.ac.uk/~gilesm/mlmc_community.html
\end{verbatim}
is also an excellent source of up-to-date information.

\subsection{Beyond European Calls and Puts}

A key step in the analysis of  section~\ref{sec:emmlmc}
was to show that the coarse and refined paths
are tightly coupled, in the sense that they
produce payoffs whose difference has small variance.
The  logic behind the analysis may be loosely summarized as
\begin{description}
\item[A] strong convergence  of Euler--Maruyama $\Rightarrow$
\item[B] coarse and refined paths close to the true path $\Rightarrow$
\item[C] coarse and refined paths close to each other
$\Rightarrow$
\item[D] coarse and refined payoffs close to each other.
\end{description}
The C $\Rightarrow$ D step appealed to the global Lipschitz
property of $h$.
This is valid
for European call and put options, where
$h(x) = \max(x-E,0)$ and
$h(x) = \max(E - x ,0)$, respectively.
However,
 the analysis must be refined for those European-style options where
$\EE[ h(X(T)) ]$ is required for functions  $h$ that violate the global Lipschitz
criterion.
We may also wish to deal with \emph{path-dependent} options where an expected value
operation is applied to a functional depending on some or all of the values
$X(t)$ for $0 \le t  \le T$.

These more exotic options include problematic classes where, for certain SDE paths,
the payoff may be very sensitive to small changes. For example, with digital options that 
expire close to the money, a
small change in the
asset path can lead to an $O(1)$ change in the payoff.
Similarly, the payoff from a barrier option is very sensitive to those paths that flirt with the
barrier.
In these cases, the logical flow above above must be adapted. Intuitively, we should be able to exploit the fact
that troublesome paths are the exception rather than the rule, and hence
C  $\Rightarrow$  D \emph{with high probability}. In some cases this allows us to recover
the computational complexity that we saw for European calls and puts.
In other cases we
must accept a slight increase in cost.

The behaviour of multilevel 
Monte Carlo for
Asian, lookback and digital options 
was considered computationally in the original 
work of Giles \cite{Giles08}.
Rigorous analysis to back up these results 
for 
 barrier, lookback and digital options 
was first given in 
 \cite{GHM09}.
Further work has been targeted at 
binary options 
\cite{AV09},
Asian options
\cite{AK14}, 
basket options \cite{GilesWS09}, 
barrier options
\cite{gdr13} and American options
\cite{BSD13}. 
The use of multilevel Monte Carlo to 
compute sensitivities with respect to problem parameters, that is, Greeks, 
was considered in 
\cite{BG12}.

\subsection{Further Developments}
It is common practice to 
combine more than one 
variance reduction technique. Given that antithetic variables 
can be effective in option valuation 
\cite{Glasserman,Hbook}, it is natural to consider embedding this 
approach within the multilevel framework.
Giles and Szpruch \cite{gs13b,GS14} have shown that this can 
be effective, particularly when Milstein is used for the numerical integration.
A conditional Monte Carlo approach has also been shown to be fruitful 
in the mutlilevel setting  
\cite{Giles07}.
In a different direction, 
Rhee and Glynn
\cite{RG2012}
have proposed an extra level of randomization 
that produces an unbiased multilevel estimator.

To go beyond the asymptotic  $O(\epsilon^{-2})$
complexity barrier it is possible to 
move to quasi Monte Carlo, 
where a low-discrepancy sequence
replaces a pseudo-random sequence.
Giles and Waterhouse 
\cite{Giles09}
have demonstrated that 
a combination of 
quasi Monte Carlo and multilevel 
can outperform each separate technique.

Finally, we note that the multilevel 
methodology has also been extended to asset models that are not driven 
purely by Brownian motion
\cite{DeHe11,XG12}.

\section{Discussion}
\label{sec:disc}
Our aim in this article was to explain in an accessible manner the 
key ideas behind the multilevel Monte Carlo method.
We focussed on the case of SDE-based financial option valuation,  
where Monte Carlo is a 
widely used tool.
At the heart of the technique is a very general and widely applicable 
philosophy---a recursive application of control variates that relies 
on tight coupling between simulations at different resolutions.
The resulting algorithm is sufficiently simple and effective that it can be 
implemented straightforwardly and used to produce noticable gains 
in computational efficiency in very general circumstances. 
However, as evidenced by the wealth of 
current research activity, there is also substantial scope for (a) refining the 
multilevel approach in order to exploit problem-specific information and (b) 
developing multilevel methods in many other stochastic simulation scenarios.
For these reasons we envisage multilevel Monte Carlo evolving into a 
cornerstone of 
computational finance.

\bigskip

\textbf{Acknowledgement}
The author is funded by a Royal Society/Wolfson Research Merit
Award and an EPSRC Digital Economy Fellowship.
He is grateful to Mike Giles for creating and placing in the public domain the code
that was used as the basis for 
	Figure~\ref{Fig.gilescode}.

\bibliographystyle{siam}
\bibliography{bookrefs,refs_ssa,Wilkie,mcmcrefs,extra_mlmc}

\end{document}